\documentclass[11pt,a4paper]{article}
\usepackage{amsfonts}

\usepackage{graphicx,setspace}
\usepackage{amsfonts,amsmath, amssymb}
\usepackage{epstopdf}
\usepackage{color,amsmath,amssymb, amsfonts,amstext,amsthm}
\usepackage{bbm}
\usepackage{mathrsfs}
\usepackage{cite}
\usepackage{amsmath}

\newcommand{\lam}{\lambda}

\newcommand{\p}{\partial}

\renewcommand{\phi}{\varphi}

\renewcommand{\a}{\alpha}

\newcommand{\R}{{\mathbb R}}

\newcommand{\st}{\theta}

\newcommand{\eps}{\varepsilon}

\newcommand{\be}{\begin{equation}}
\newcommand{\ee}{\end{equation}}
\newcommand\bes{\begin{eqnarray}}
\newcommand\ees{\end{eqnarray}}
\newcommand{\bess}{\begin{eqnarray*}}
\newcommand{\eess}{\end{eqnarray*}}

\newcommand{\dx}{{\rm d}x}
\newcommand{\dy}{{\rm d}y}

\newcommand{\befig}{\begin{figure}}
\newcommand{\enfig}{\end{figure}}

\newcommand{\bear}{\begin{eqnarray}}
\newcommand{\enar}{\end{eqnarray}}

\newcommand{\bearn}{\begin{eqnarray*}}
\newcommand{\enarn}{\end{eqnarray*}}
 
{\setlength\arraycolsep{2pt}}

\title{Mean exit time  for stochastic dynamical systems driven by tempered stable L\'evy fluctuations
}
\author{Yanjie Zhang$^1$,  Xiao Wang$^2$ \footnote{Corresponding author}, and Jinqiao Duan$^3$  \\
\\
\ \\
   {\small \it $^1$ School of Mathematics, South China University of Technology }\\
  {\small \it Guangzhou 510000,  China }\\
  {\small \tt email:zhangyj18@scut.edu.cn}\\
  {\small \it $^2$ School of Mathematics and Statistics, Henan University }\\
  {\small \it Kaifeng 475001, China  }\\
  {\small \tt email: xwang@vip.henu.edu.cn}\\
  {\small \it $^3$ Department of Applied Mathematics, Illinois Institute of Technology }\\
  {\small \it Chicago, IL 60616, USA }\\
   {\small \tt email:duan@iit.edu }
}
\begin{document}
\maketitle

\begin{abstract}
 We use the mean exit time to quantify macroscopic dynamical behaviors of stochastic dynamical systems driven by tempered L\'evy fluctuations,  which are solutions of nonlocal elliptic equations. Firstly, we construct a new numerical scheme to compute and solve the mean exit time associated with the one dimensional stochastic system. Secondly, we extend the analytical and numerical results to two dimensional case: horizontal-vertical and isotropic case. Finally, we verify the effectiveness of the presented schemes  with numerical experiments  in  several examples.


\end{abstract}

\baselineskip=15pt

\section{Introduction}\label{intro}

 Because of the boundedness of the physical space, the extremely heavy tails of these models are not realistic for most real-world applications. This has led researchers to use models that are similar to stable distributions in some central region, but with lighter tails. Tempered stable distributions are a  class of models that capture this type of behavior, which describe the trapped dynamics, widely appearing in nature \cite{Koponent95, Rosinski07}.

The mean exit time (MET) is an important tool to quantify macroscopic dynamical behaviors of a stochastic system, as it describes the expected time of a particle initially inside a bounded domain until the particle first exits the domain. Deng et al. studied the  mean exit time for the anomalous processes having the tempered L\'evy stable waiting times in the theory \cite{WDX, Deng}. Motivated the previous work, in this letter, we construct new numerical schemes to compute and solve the mean exit time associated with these one and two dimensional stochastic systems. Furthermore, we verify the effectiveness of the presented schemes  with numerical experiments  in  several examples.

\section{MET for one-dimensional case }
\label{1dim}
Consider the following one dimensional stochastic dynamical system
\begin{equation}
\label{SDE011}
\mathrm{d}X_t=f(X_t)\mathrm{d}t+ \mathrm{d}L_t,
\end{equation}
where $f$ is a drift term (vector field), and $L_t$ is a tempered stable L\'evy process  with triplet $(0,d, \kappa\nu)$. i.e., zero linear coefficient, diffusion coefficient $d\geq 0$, L\'evy measure $\kappa\nu(dy)$ and $\kappa$ is a nonnegative parameter. The  jump measure $\nu$ for one dimensional tempered L\'evy process is obtained by multiplying the $\a$-stable L\'evy measure $\nu_{\alpha}(dy)$ by an exponential decaying function, i.e.,

\bear  \label{LevyMesure}
   \nu(dy) = \nu_{\a}(dy) \left(1_{\{y>0\}}e^{-\lam_1y}+1_{\{y<0\}}e^{\lam_2y}\right)
            =  \left[ \frac{C_{\alpha,\lambda_1}}{e^{\lam_1y}y^{1+\a}}1_{\{y>0\}}+\frac{C_{\alpha,\lambda_2}}{e^{-\lam_2y}(-y)^{1+\a}}1_{\{y<0\}}\right]dy,
\enar
where $C_{\alpha, \lambda_{i}} (i=1,2)$ is a positive constant,  $\a \in(0,1) \bigcup (1,2)$ is called the stable index, and $ \lambda_i$ is the positive  tempering  parameter. Here we consider `symmetric'  tempered L\'evy process, i.e., $\lam_1=\lam_2=\lambda$, then   $C_{\alpha, \lambda_i}=C_\a=\frac{1}{2|\Gamma(-\a)|}$ (see \cite{Deng}).
\par
 The mean exit time for the solution orbit $X_t$ in Eq.~(\ref{SDE011}) starting at $x$ from a bounded domain $D$ is defined as
\begin{equation}
 \tau_x{(\omega)}:=\inf \{t \geq 0: X_{t}(\omega , x) \notin  D, X_0=x\},~~
u(x):=\mathbb{E}[\tau_x(\omega)],
\end{equation}
which satisfies the following integro-differential equation
\begin{equation}
  \label{MET}
       \mathscr{L}u = -1, ~x \in D,~~u(x) = 0,  ~~x \in D^{c},
    \end{equation}
where
\begin{equation}
\label{MET0}
 \mathscr{L}u =f(x)u_x+\frac{d}{2}u_{xx}  + \eps\int_{\mathbb{R}\setminus \{0\}}
   \left[u(x+y)-u(x)+1_{\{|y|< 1\}}(y)y u_x\right]\nu(\mathrm{d}y).
\end{equation}
In the following, we will construct a new numerical scheme to compute the MET for one dimensional stochastic dynamical system with a scalar tempered L\'evy fluctuation.

\subsection{Numerical schemes}

\par
Introduce the following function,
\begin{equation}
 \int_s^{\infty} x^{-\varrho} e^{- x}\dx = s^{-\frac{\varrho}{2}}e^{-\frac{s}{2}}W_{-\frac{\varrho}{2},\frac{1-\varrho}{2}}(s), ~~for~~ s>0,
\end{equation}
where $W $ is the Whittaker W function.

Assume the spatial domain $D=(-1,1)$, in the sense of the principal value,  the integral 
  $\int_{\mathbb{R}\setminus \{0\}} 1_{\{|y|<1\}}(y)y u_x \nu(\mathrm{d}y)$ vanishes,  
  then the integral term  of equation \eqref{MET0} becomes
\begin{equation}
\label{FD}
\begin{aligned}
\int_{\mathbb{R}\setminus \{0\}}\left[u(x+y)-u(x)\right]\nu(\mathrm{d}y)&= C_{\a}\int_{-\infty}^{-1-x} \frac{u(x+y)-u(x)}{e^{\lam |y|}|y|^{1+\a}} \dy +C_{\a}\int_{-1-x}^{1-x} \frac{u(x+y)-u(x)}{e^{\lam |y|}|y|^{1+\a}} \dy\\
&+ C_{\a}\int_{1-x}^{\infty} \frac{u(x+y)-u(x)}{e^{\lam |y|}|y|^{1+\a}} \dy   \\
&= -C_{\a} u(x) \left[W_1(x)+W_2(x) \right]+C_{\a} \int_{-1-x}^{1-x} \frac{u(x+y)-u(x)}{e^{\lam |y|}|y|^{1+\a}} \dy,
\end{aligned}
\end{equation}
where
\bess
    W_1(x)&= \lam^{\frac{\a-1}{2}}(1+x)^{-\frac{\a+1}{2}} e^{-\frac{\lam (1+x)}{2}} W_{-\frac{1+\a}{2},-\frac{\a}{2}}(\lam(1+x)), \\
      W_2(x)&= \lam^{\frac{\a-1}{2}}(1-x)^{-\frac{\a+1}{2}} e^{-\frac{\lam (1-x)}{2}} W_{-\frac{1+\a}{2},-\frac{\a}{2}}(\lam(1-x)).
\eess
\par
For the singular integral term of equation \eqref{FD}, we take $\delta=\min\{1-x,1+x\}$,  using a modified trapezoidal rule for the singular term, then we have
\begin{equation}
\begin{aligned}
      \int_{-1-x}^{1-x} \frac{u(x+y)-u(x)}{e^{\lam |y|}|y|^{1+\a}} \dy &=\mathbb{P.V.} \int_{-1-x}^{1-x} \frac{u(x+y)-u(x)-1_{\{|y|<\delta\}}yu'(x)}{e^{\lam |y|}|y|^{2}} \big|y\big|^{1-\a}\dy\\
          &= \int_{0}^{1-x} g(y)y^{1-\a}\dy +\int_{0}^{1+x} \tilde{g}(y)y^{1-\a}\dy \\
          &= h \sum\limits_{j=1}^{J_{1}}\!{'} G(y_j) - \zeta(\a-1)g(0)h^{2-\a} - \zeta(\a-2)g'(0)h^{3-\a} +O(h^2)   \\
          &~~~+ h \sum\limits_{j=1}^{J_{2}}\!{'} \tilde{G}(y_j) - \zeta(\a-1)\tilde{g}(0)h^{2-\a}-\zeta(\a-2)\tilde{g}'(0)h^{3-\a} +O(h^2),
\end{aligned}
\end{equation}
where $g(y)=\frac{u(x+y)-u(x)-1_{\{|y|<\delta\}}yu'(x)}{e^{\lam |y|}|y|^{2}}$,
$\tilde{g}(y)=g(-y)$, $G(y) =g(y)|y|^{1-\a}$, $\tilde{G}(y)=G(-y)$, $J_{1}$ and $J_{2}$ are the index corresponding to $1-x$ and $1+x$, respectively.
Moreover, $h\cdot J_{1}=1-x$ and $h\cdot J_{2}=1+x$. The summation symbol $\displaystyle{\sum\limits\!{'}}$ means the term of upper index is multiplied by $\frac{1}{2}$, $\zeta$ is the Riemann zeta function, $g(0) = \tilde{g}(0)=\frac{u''(x)}{2}$,
$g'(0)= \frac{u'''(x)}{6}-\lam g(0)$,
$\tilde{g}'(0) = -\frac{u'''(x)}{6}+\lam g(0)$.

Let us divide the interval $[-2,2]$ into $4J$ subintervals and define $x_j=jh$ for $-2J\leq j\leq 2J$ integer, where $h=\frac{1}{J}$. Using central difference numerical scheme for the first and two derivatives and modifying the ``punched-hole" trapezoidal rule in the nonlocal term, we get the discretization scheme of (\ref{MET}), i.e.,
\begin{equation}
\label{eq:Discrete}
\begin{aligned}
&C_h\frac{U_{j+1}-2U_j+U_{j-1}}{h^2} - f(x_j) \left(\frac{U_{j+1}-U_{j-1}}{2h}\right)
     -\kappa C_{\a} \left[W_1(x_j)+W_2(x_j)\right]U_j\\
        &+ \kappa C_{\a} h \sum\limits\!{''}_{k=-J-j, k\neq 0}^{J-j}\frac{U_{j+k}-U_j}{e^{\lam |x_k|}|x_k|^{1+\a}} =-1,
\end{aligned}
\end{equation}
where $\displaystyle{\sum\limits\!{''}}$ means that the quantities corresponding to the two end summation indices are
multiplied by 1/2 and $C_h = \frac{d}{2} -\eps C_{\a} \zeta(\a-1)h^{2-\a}$.

We can rewrite the summation terms of Eq.~(\ref{eq:Discrete}) as multiplication form of matrix-vector $R\mathbf{U}$, where
 $R$ is a $(2J-1)\times(2J-1)$ matrix. Moreover, the matrix $R$ can be decomposed as
 \begin{equation}
 R = T_R + D_R,
 \end{equation}
 where $T_R$ is a Toeplitz matrix , i.e.,

\bess
T_R =
\begin{pmatrix}
	0 & \frac{\widetilde{C}}{e^{\lam h}h^{1+\a}} & \frac{\widetilde{C}}{e^{2\lam h}(2h)^{1+\a}} &\cdots & \cdots & \frac{\widetilde{C}}{e^{(2J-2)\lam h}[(2J-2)h]^{1+\a}}\\
	\frac{\widetilde{C}}{e^{\lam h}h^{1+\a}} & 0 &\frac{\widetilde{C}}{e^{\lam h}h^{1+\a}} &\cdots& \cdots  & \frac{\widetilde{C}}{e^{(2J-3)\lam h}[(2J-3)h]^{1+\a}}\\
	\frac{\widetilde{C}}{e^{2\lam h}(2h)^{1+\a}} & \frac{\widetilde{C}}{e^{\lam h}h^{1+\a}} & 0  & \ddots & \cdots & \frac{\widetilde{C}}{e^{(2J-4)\lam h}[(2J-4)h]^{1+\a}} \\
	\vdots & \vdots & \ddots & \ddots& \ddots  & \vdots\\
\frac{\widetilde{C}}{e^{(2J-3)\lam h}[(2J-3)h]^{1+\a}} &  \frac{\widetilde{C}}{e^{\lam h}[(2J-4)h]^{1+\a}} & \cdots & \ddots & 0 & \frac{\widetilde{C}}{e^{\lam h}h^{1+\a}} \\
	\frac{\widetilde{C}}{e^{(2J-2)\lam h}[(2J-2)h]^{1+\a}} &  \frac{\widetilde{C}}{e^{(2J-3)\lam h}[(2J-3)h]^{1+\a}} & \frac{\widetilde{C}}{e^{(2J-4)\lam h}[(2J-4)h]^{1+\a}} & \cdots  & \frac{\widetilde{C}}{e^{\lam h}h^{1+\a}} & 0
\end{pmatrix}
\eess

and $D_R$ is a tridiagonal one, i.e.,
\bess
D_{R} =
\begin{pmatrix}
  a_{1-J} & 0 &         &       &      \\
 0 & a_{2-J}  & 0  &       &      \\
        & \ddots & \ddots  & \ddots & \\
	&        & 0 & a_{J-2} & 0\\
	&        &           &0 & a_{J-1}

\end{pmatrix}
\eess
with
$$\widetilde{C} = \kappa C_{\a} h, ~ a_l = - \eps C_{\a} h \sum\limits\!{''}_{k=-J-l, k\neq 0}^{J-l}\frac{1}{e^{\lam |x_k|}|x_k|^{1+\a}},~
      l=1-J, 2-J \cdots, J-1.$$

\subsection{Numerical experiments}

\subsubsection{ Verification}

Taking $u(x)=(1-x^2)_+$ ( i.e., $u(x)=1-x^2$ for $x \in (-1,1)$, otherwise, $u(x)=0$)
and $\lam_1=\lam_2=\lambda, f=d=0,\eps=1$ into the right-hand side (RHS) of Eq.~(\ref{MET0}), we have
\bear
    RHS&=&C_{\alpha}\int_{-1-x}^{1-x} \frac{-2xy-y^2}{e^{\lam |y|}|y|^{1+\a}} \dy
        -C_{\alpha}u(x)\left[\int_{1-x}^{\infty} \frac{\dy}{e^{\lam y}y^{1+\a}} +\int_{1+x}^{\infty} \frac{\dy}{e^{\lam y}y^{1+\a}} \right]  \nonumber\\
            &=& 2C_{\alpha}x \lam^{\frac{\a}{2}-1}[(1-x)^{-\frac{\a}{2}}e^{-\frac{\lam (1-x)}{2}}W_{-\frac{\a}{2},\frac{1-\a}{2}}(\lam (1-x))
             -(1+x)^{-\frac{\a}{2}}e^{-\frac{\lam (1+x)}{2}}W_{-\frac{\a}{2},\frac{1-\a}{2}}(\lam (1+x))]  \nonumber\\
            &&-C_{\alpha}\Gamma(2-\a)\lam^{\a-2}\left[P(2-\a,\lam(1-x))+P(2-\a,\lam(1+x))] -C_{\alpha}(1-x^2)[W_1(x)+W_2(x)\right],   \nonumber
\enar
where  $P(a,x)= \frac{\int_0^{x} e^{-y}y^{a-1} \dy}{\Gamma(a)}(a \geq 0) $ is the incomplete Gamma function, and
$Q(a,x)=1-P(a,x)$ is the `upper' incomplete Gamma function.

\begin{figure}
\begin{minipage}[t]{0.4\linewidth}
\centering
\includegraphics[width=3.5in]{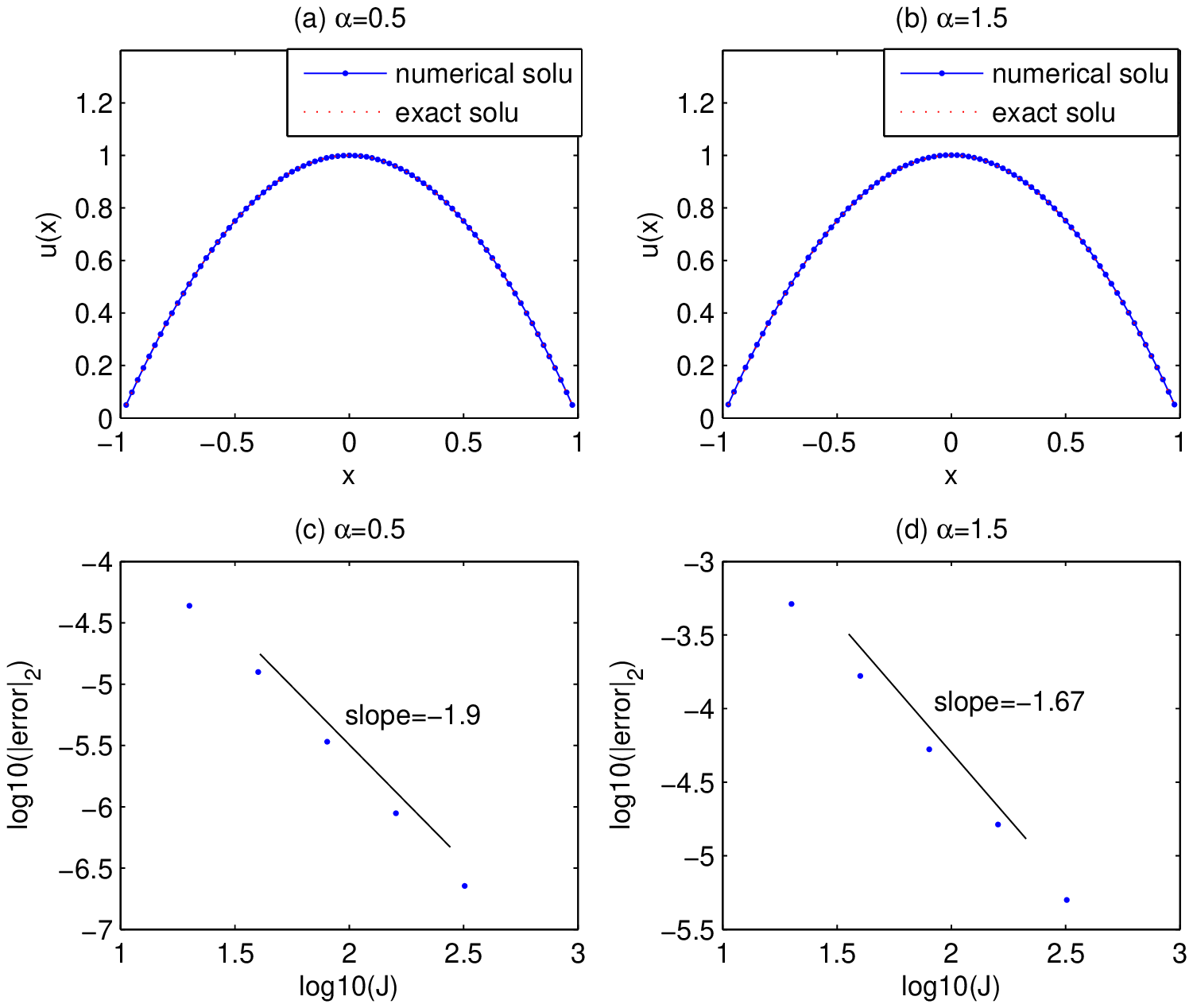}
\caption{(a)   Comparison between numerical solution  and exact
solution $u(x)=(1-x^2)_+$ for $\a=0.5, \lam=0.01$; (b)the same as (a) except $\a=1.5$;
(c) the error between numerical solution and exact solution for $\a=0.5$.
(d) the same as (c) except $\a=1.5$.  }
\label{ComExpOrd}
\end{minipage}%
 \hspace{0.8in}
\begin{minipage}[t]{0.3\linewidth}
\centering
\includegraphics[width=2.8in]{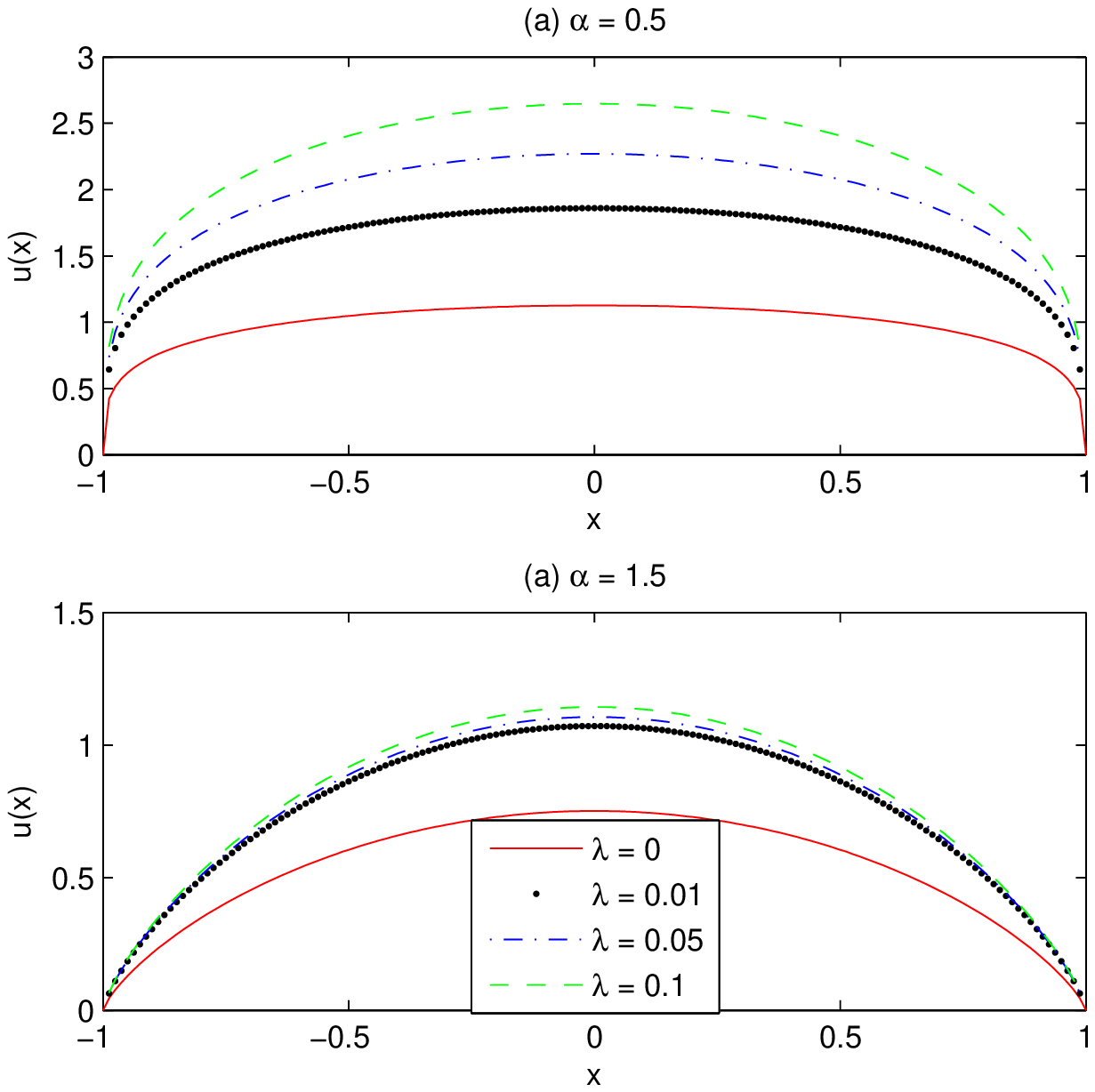}
\caption{The solutions of mean exit time $u(x)$ of Eq.~(\ref{MET}) for different $\lam$. (a) $\a=0.5$; (b)$\a=1.5$. }
\label{MET1DExp}
\end{minipage}
\end{figure}

Take the exact solution $u(x)=(1-x^2)_+$ of constructed equation to verify our numerical method and
compute the {\bf{convergence orders}}. Fig.~\ref{ComExpOrd} shows the errors between the numerical and exact solutions
with $\lam=0.01, ~f=d=0,~\eps=1$ and different $\alpha$. Fig.~\ref{ComExpOrd}(a) and (b) show our numerical solution almost agree with the exact
solutions for different $\a$($\a=0.5,~ \a=1.5$). The numerical
convergence order is equal to 2. To verify it, we plot $\log_{10}(|error|_2)$ against $\log_{10}(J)$ with different
resolutions $J=20,40,80,160,320$ in Fig.~\ref{ComExpOrd}(c) and Fig.~\ref{ComExpOrd}(d), where $|error|_2$ represents the 2-norm errors.  This above results imply that the errors almost reach our order expected from the above analysis.

\befig[h]
\begin{center}
\includegraphics*[width=0.6\linewidth]{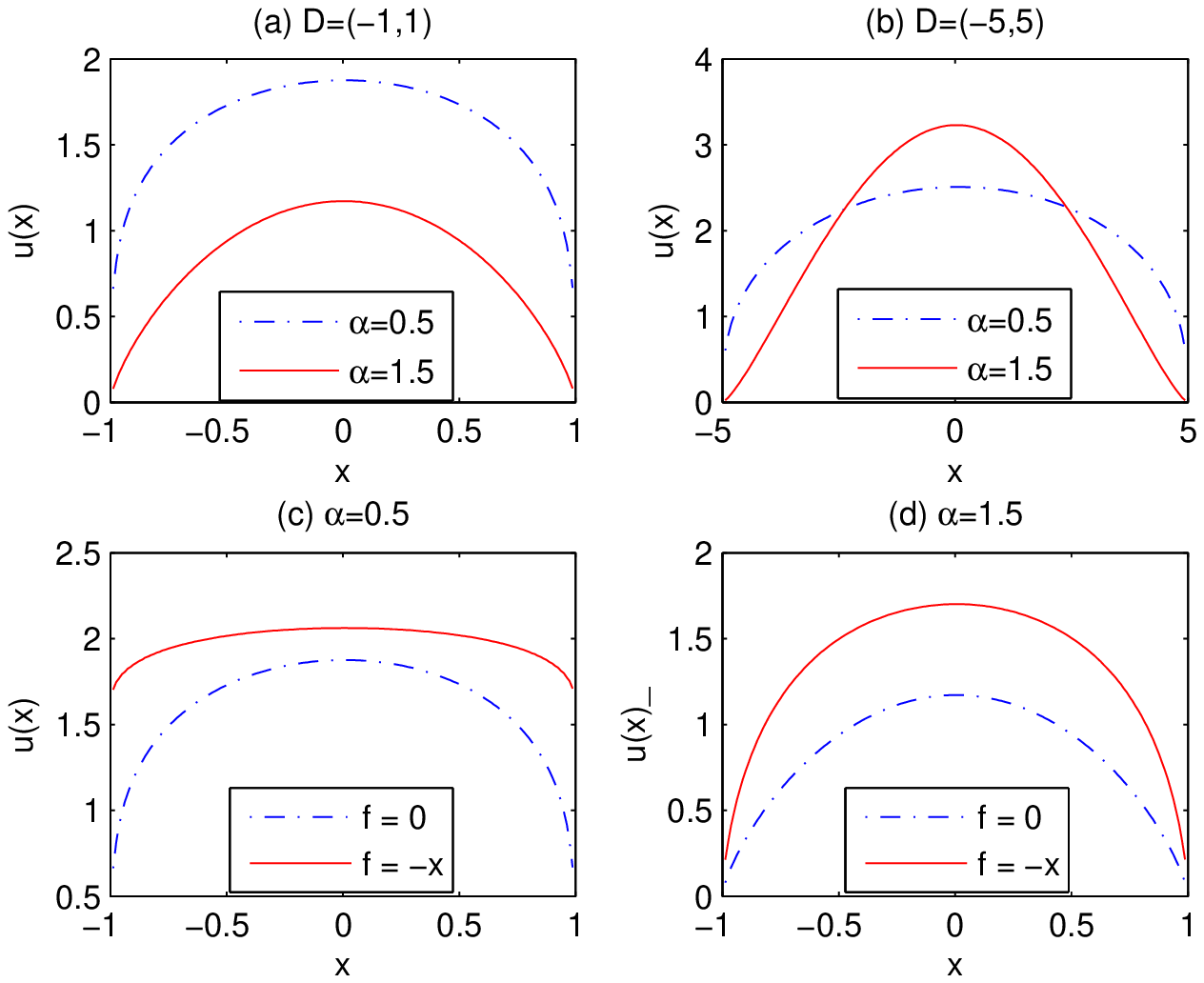}
\caption{The effect of domain $D$ and drift term $f$ on MET $u(x)$ of Eq.~(\ref{MET}) for $\lam =  0.01$ and
$d=0, \eps=1$. (a) the domain $D=(-1,1)$ for different $\a=0.5,1.5$ with $f=0$; (b) the same as (a) except
$D=(-5,5)$;(c) $\a=0.5, D=(-1,1)$ for different drift term $f$; (d) the same as (c) except $\a=1.5$ . }
\label{DifD_f}
\end{center}
\enfig

\subsubsection{ Effect of parameters}
Here we consider the effect of tempering parameter $\lam$ for MET. Fig.~\ref{MET1DExp} shows the numerical solution of MET for different $\lam $~($\lam = 0, 0.01, 0.05, 0.1$) and $\a$~($\a=0.5, 1.5$) with $f=d=0, ~\eps=1,~ D=(-1,1)$. For $\lam=0$,  we use the method in reference \cite{Tinge} for comparison.
For $\a=0.5$~(see Fig.~\ref{MET1DExp}(a)),
the `particle' takes more time to exit as  $\lam$ becomes larger, which  agrees
with our intuition, i.e., the L\'evy measure becomes smaller as the tempering parameter $\lam$ becomes larger, then the jump
intensity is smaller and the `particle' is harder to exit the domain. Fig.~\ref{MET1DExp}(b) shows the similar
results, but the effect of tempering parameter is small for $\a=1.5$. It is also interesting to point out the effect of domain $D$ and drift term $f$ for MET. When the other parameters
are fixed, we find that the `particle' will take more time to exit the domain as the domain becomes larger
in Fig.~\ref{DifD_f} (a) and (b). For $D=(-5,5)$, we find that the MET increases when the parameter $\alpha$ increases near the origin. However,
for $D=(-1,1)$, the MET decreases when the parameter $\alpha$ increases near the origin.
In Fig.~\ref{DifD_f} (c) and (d), the `particle' is harder to exit the domain, because the drift term `$f(x)=-x$' drives it toward the origin.

\section{MET for two-dimensional case}
Consider the following two dimensional stochastic dynamical system
\begin{equation}
\label{SDE01}
\mathrm{d}X_t=f(X_t)\mathrm{d}t+ \mathrm{d}L_t,
\end{equation}
where $f$ is a vector field, and $L_t$ is a tempered stable L\'evy process  with triplet $(0,\mathbf{d}, \kappa\nu)$, $\mathbf{d}$ is a symmetric non-negative definite matrix, the jump measure $\nu$ is the following two cases: horizontal-vertical case and isotropic case, i.e.,
$
  {\nu}(\mathrm{d}y) = \frac{C_1}{e^{\lam_1 y_1}|y_1|^{1+\a_1}}\delta(y_2)\mathrm{d}y_1\mathrm{d}y_2+\frac{C_2}{e^{\lam_2 y_2}|y_2|^{1+\a_2}}\delta(y_1)\mathrm{d}y_1\mathrm{d}y_2$ and $\nu(dy) = \frac{\widetilde{C}_{\a} dy}{ e^{\lam |y|}|y|^{\alpha+2}}$
with $C_1=\frac{1}{2|\Gamma(-\a_1)|}$, $C_2=\frac{1}{2|\Gamma(-\a_2)|}$ and $\widetilde{C}_{\a} = \frac{1}{2\pi|\Gamma(-\a)|}$.

The usual exponentially tempered L\'evy measure $\nu$ is expressed as (see \cite{JR})
\bear \label{LevyMesureGen}
   { \nu} (B) = \int_{S_2}\Gamma({\rm d}\theta)\int_0^{\infty}1_B(r\theta)e^{-r}{r^{-1-\a}}{\rm d}r, \; \forall B \in  \mathscr{B}(\mathbb{{R}}^2)
\enar
with  $S_2=\{x:|x|=1\}$   the unit circle in $\mathbb{R}^2$,  and $\Gamma$ is the finite measure on this unit circle.

The generator for (\ref{SDE01}) is
\begin{equation}
\begin{aligned}
\label{METexpLey2D}
\widetilde{\mathscr{L}}u(x) &=f^{i} (\p_i u)(x)+\frac12 d^{ij}(\p_i \p_ju)(x)+ \kappa \int_{\mathbb{R}^2\setminus \{\mathbf{0}\}}
       \left[u(x+y)-u(x)+1_{B_h(\mathbf{0})}y^{i} (\p_{i}u)(x)\right]\nu(\mathrm{d}y),
\end{aligned}
\end{equation}
where $B_h(\mathbf{0}) = \{x:|x| \leq h \ll 1\}$.

\subsection{MET for the horizontal-vertical case}

When   the components of  the  tempered L\'evy process $L_t$ are   independent,   the particles (or solutions)  spread in either horizontal or vertical direction \cite{Deng}. The finite measure $\Gamma$ in (\ref{LevyMesureGen}) concentrates on the points of intersection of unit circle  $S_2$ and axes. The MET satisfies the following integro-differential equation
\begin{equation}
 \label{twod}
    \widetilde{ \mathscr{L}}u  = -1, ~~ x  \in D,~~
               u(x) = 0,   ~~x \in D^{c}.
\end{equation}

\subsubsection{Numerical methods}

Here we take $\a_1=\a_2=\alpha, \lam_1=\lam_2=\lambda$, $C_1=C_2=C_\a$ and the square domain $D=(-1,1)^2$. Set $x=(x_1,x_2)\in \mathbb{R}^{2}$ and $y=(y_1,y_2)\in \mathbb{R}^{2}$. The integral terms in (\ref{twod}) can be divided into two parts, i.e.,
\bear\label{FD2D}
   &&\int_{\mathbb{R}^2\setminus \{\mathbf{0}\}}\left[u(x_1+y_1,x_2+y_2)-u(x_1,x_2)\right]{\nu}(\mathrm{d}y) \nonumber  \\
     &=& -C_\a u(x_1,x_2) [W_1(x_1)+W_2(x_1) ]+
             C_\a\int_{-1-x_1}^{1-x_1} \frac{\left[u(x_1+y_1,x_2)-u(x_1,x_2)\right]}{e^{\lam |y_1|}|y_1|^{1+\a}}\dy_1     \nonumber  \\
            &&-C_\a u(x_1,x_2) [W_1(x_2)+W_2(x_2) ]+
             C_\a \int_{-1-x_2}^{1-x_2} \frac{\left[u(x_1,x_2+y_2)-u(x_1,x_2)\right]}{e^{\lam |y_2|}|y_2|^{1+\a}}\dy_2.
\enar

Similarly, we use the modified trapezoidal rule  for the integral terms in (\ref{FD2D}) to get
\bess
   \int_{-1-x_i}^{1-x_i}\tilde{G}(y_i)\dy_i
     =  h \sum\limits\!{'}_{k=-J-j, k\neq0}^{J-j} \tilde{G}(y_{i_k}) -\zeta(\a-1)h^{2-\a} u_{x_ix_i}(x_1,x_2)+O(h^2), ~~i=1,2,
\eess
where
\begin{equation}
\tilde{G}(y_1)=\frac{u(x_1+y_1,x_2)-u(x_1,x_2)}{e^{\lam |y_1|}|y_1|^{1+\a}},  ~~
\tilde{G}(y_2)=\frac{u(x_1,x_2+y_2)-u(x_1,x_2)}{e^{\lam |y_2|}|y_2|^{1+\a}}.
\end{equation}

\subsubsection{Numerical experiments}
\par
 Here we fix the factors $f^i=0, d^{ij}=0, D=(-1,1)^{2}, \eps=1$. Fig.~\ref{Met2D1} displays the MET for two-dimensional horizontal-vertical case with different $\lam$ and $\a$. We find that the MET increases as the parameter $\lam$ increases. However, it decays faster for $\a=1.5$ than $\a=0.5$ near the boundary.

\subsection{MET for the isotropic case}
When the particles spread uniformly in all directions, this case is called the isotropic L\'evy process. Here we assume the process is radially symmetric and the domain $D=\{x \in \mathbb{R}^{2}:|x|<1\}$, then we have $u(x)= u(r)$, where $x=(x_1,x_2)$ and $r = |x| = \sqrt{x_1^2+x_2^2}$.

Set $\frac{d^{ij}}{2} = d(r)\mathbf{I}$ and $f^i=f(r)\frac{x_i}{r}$, $i=1,2$,
where $f(\cdot)$ and $d(\cdot)$ are smooth scalar functions, then the MET satisfies the following integro-differential equation
\begin{align}
\label{eq:equivform}
    &f(r) u'(r) + d(r)\left[u''(r)+\frac{u'(r)}{r}\right]+\kappa \widetilde{C}_{\a}\int_{\R^2 \setminus\{\mathbf{0}\}}\frac{u(x+y)-u(x)-1_{B_{h}(\mathbf{0})}y^{i} (\partial_i u)(x)}{e^{\lam |y|}|y|^{\alpha+2}}{\rm d}y =-1.
\end{align}

\subsubsection{Numerical methods}

\par
For the radially symmetric case, we only consider the solution $u(x)$ on the positive $x_1$-axis. For simplicity, we denote $x=(r,0)$ for $r \ge 0$.
By taking $0<h\ll 1$, the singular integral term in Eq.(\ref{eq:equivform})becomes
\begin{align} \label{split}
     &\int_{\mathbb{R}^2\setminus \{\mathbf{0}\}
    }\frac{u(x+y)-u(x)-1_{B_{h}(\mathbf{0})}y^i (\partial_i u)(x) u(x)}{e^{\lam |y|}|y|^{\alpha+2}} {\rm d}y   \nonumber \\
   &=\int_{\mathbb{R}^2\setminus B_{h}(\mathbf{0})}\frac{u(x+y)-u(x)}{e^{\lam |y|}|y|^{\alpha+2}} {\rm d}y +
     \int_{B_{h}(\mathbf{0}) \setminus \{\mathbf{0}\}  }\frac{u(x+y)-u(x)-y^i (\partial_i u)(x)}{e^{\lam |y|}|y|^{\alpha+2}} {\rm d}y\\
   &= 2\int_{(0,1)\setminus(r-h,r+h)} s[u(s)-u(r)]F_{\lam}^1(s,r)  {\rm d}s
    +2\int_{(r-h,r+h)} s[u(s)-u(r)] F_{\lam}^2(s,r) {\rm d}s  \\
        &~~~~-2\lam^{\frac{\a-1}{2}} u(r)\int_0^\pi \tilde{r}^{-\frac{\a+1}{2}} e^{-\frac{\lam \tilde{r}}{2}} W_{-\frac{1+\a}{2},-\frac{\a}{2}}(\lam\tilde{r})  {\rm d} \st
        +C_0\left[u''(r)+\frac{u'(r)}{r}\right]+ \mathcal{O}(h^{4-\a})
        \end{align}
where $P(a,x) $ is the incomplete Gamma function, $C_0=\pi \lam^{2-\a} \Gamma(2-\a) P(2-\a, \lam h)$,  $\tilde{r}=\sqrt{1-r^2\sin^2{\st}}-r\cos{\st}$, and
\begin{equation}
\begin{aligned}
F_{\lam}^1(s,r) &=  \int_{0}^\pi e^{-\lam \sqrt{s^2+r^2-2sr\cos{\st}}}[s^2+r^2-2sr\cos{\st}]^{-\frac{\a+2}{2}} {\rm d}\st,\\
F_{\lam}^2(s,r) &=  \int_{\gamma}^\pi e^{-\lam \sqrt{s^2+r^2-2sr\cos{\st}}}[s^2+r^2-2sr\cos{\st}]^{-\frac{\a+2}{2}} {\rm d}\st.
\end{aligned}
\end{equation}

For $r\neq 0$, the integro-differential equation \eqref{eq:equivform} can be rewritten as
\begin{equation}
\label{wholeI}
\begin{aligned}
&f(r) u'(r) + d(r)\left[u''(r)+\frac{u'(r)}{r}\right]+2\kappa \widetilde{C}_{\a}\int_{(0,1)\setminus(r-h,r+h)} s[u(s)-u(r)]F_{\lam}^1(s,r)  {\rm d}s\\
 &+2\kappa \widetilde{C}_{\a}\int_{(r-h,r+h)} s[u(s)-u(r)] F_{\lam}^2(s,r) {\rm d}s
        -2\kappa \widetilde{C}_{\a}\lam^{\frac{\a-1}{2}} u(r)\int_0^\pi \tilde{r}^{-\frac{\a+1}{2}} e^{-\frac{\lam \tilde{r}}{2}} W_{-\frac{1+\a}{2},-\frac{\a}{2}}(\lam\tilde{r})  {\rm d} \st\\
        &+\kappa C_0 \widetilde{C}_{\a}\left[u''(r)+\frac{u'(r)}{r}\right]+ \mathcal{O}(h^{4-\a})=-1, .
\end{aligned}
\end{equation}

 For $r=0$, we have
\begin{equation}
\begin{aligned}
  &f(0)u'(0)+\left(d(0)+\kappa \widetilde{C}_{\a}C_0\right)\left[\frac{\p^2u}{\p x_1^2}+\frac{\p^2u}{\p x_2^2}\right]\bigg\lvert_{x=0}
      + 2\pi\kappa\widetilde{C}_{\a}\int_{h}^1 \frac{u(r)-u(0)}{e^{r\lam}r^{\a+1}}{\rm d}r -2\pi \kappa \widetilde{C}_{\a} W_1(0)u(0)\\
      &=-1.
\end{aligned}
\end{equation}

\subsubsection{Numerical experiments}

\begin{figure}
\begin{minipage}[t]{0.5\linewidth}
\centering
\includegraphics[width=3.5in]{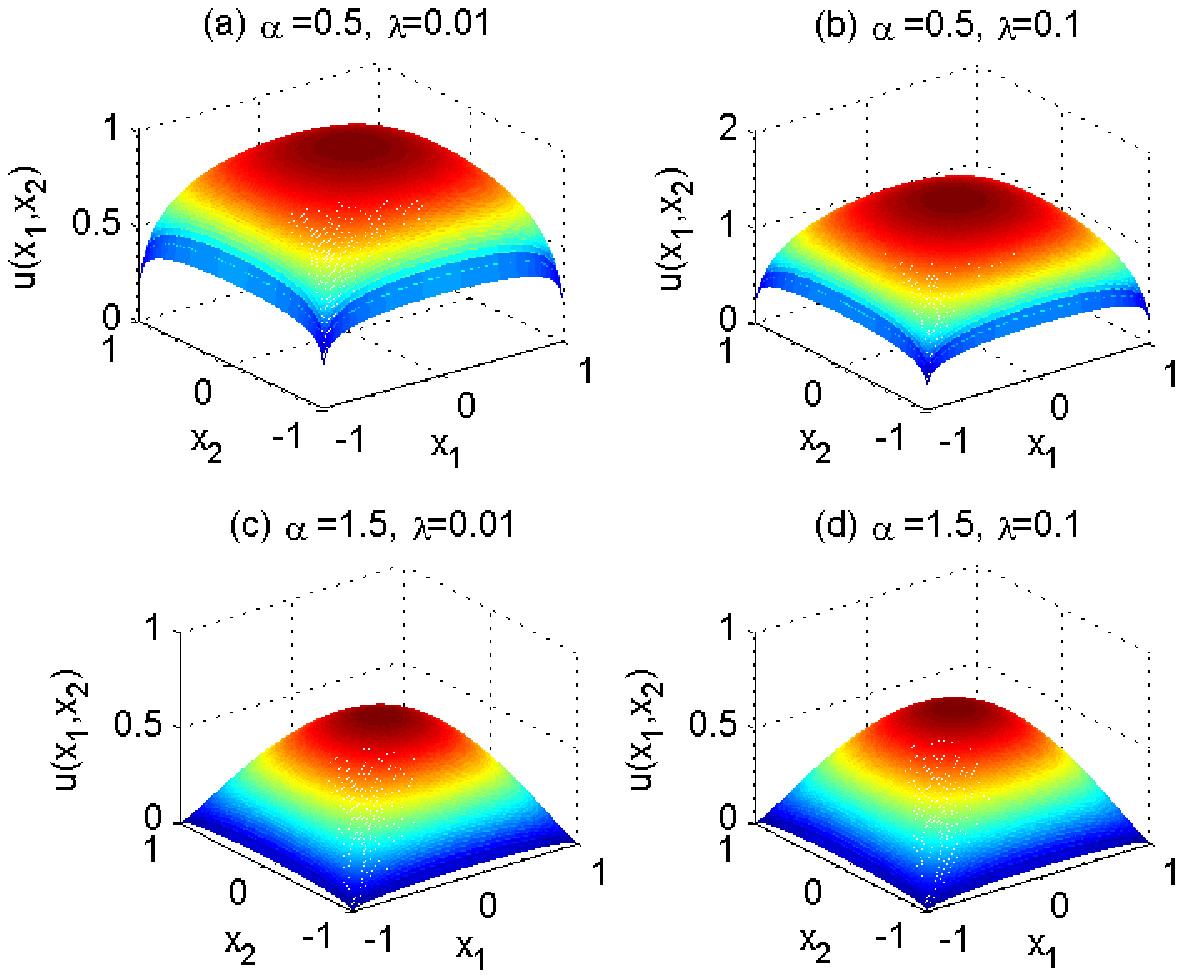}
\caption{MET for the horizontal-vertical case with different $\lam$ ($\lam = 0.01, 0.1$) and $\a$ ($\a=0.5, 1.5$).  }
\label{Met2D1}
\end{minipage}%
 \hspace{0.5in}
\begin{minipage}[t]{0.4\linewidth}
\centering
\includegraphics[width=3.1in]{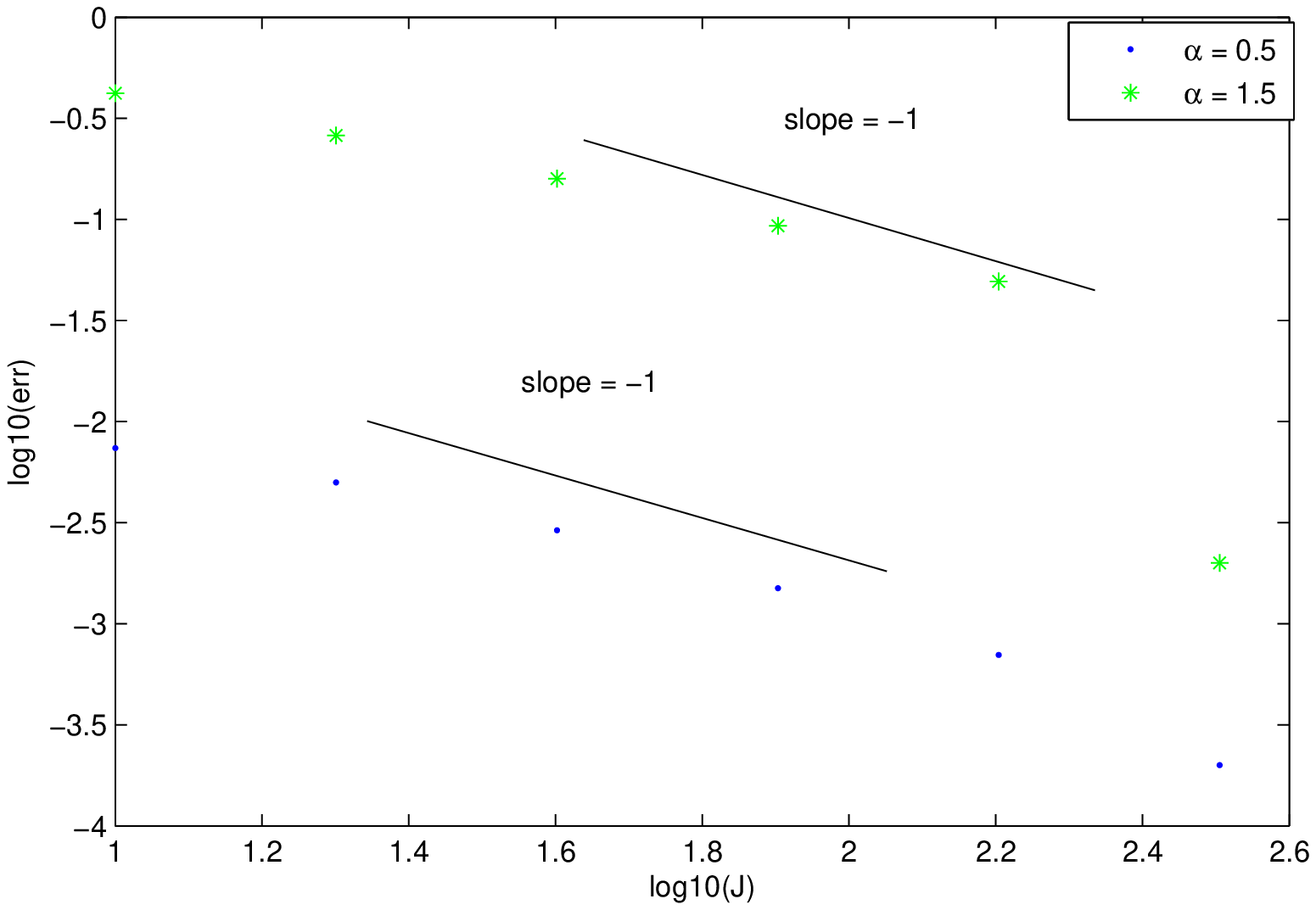}
\caption{The order of MET for the isotropic case with $\lam=0.01$ and $\a$.  }
\label{order2D}
\end{minipage}%
\end{figure}

\befig[h]
\begin{center}
\includegraphics*[width=0.7\linewidth]{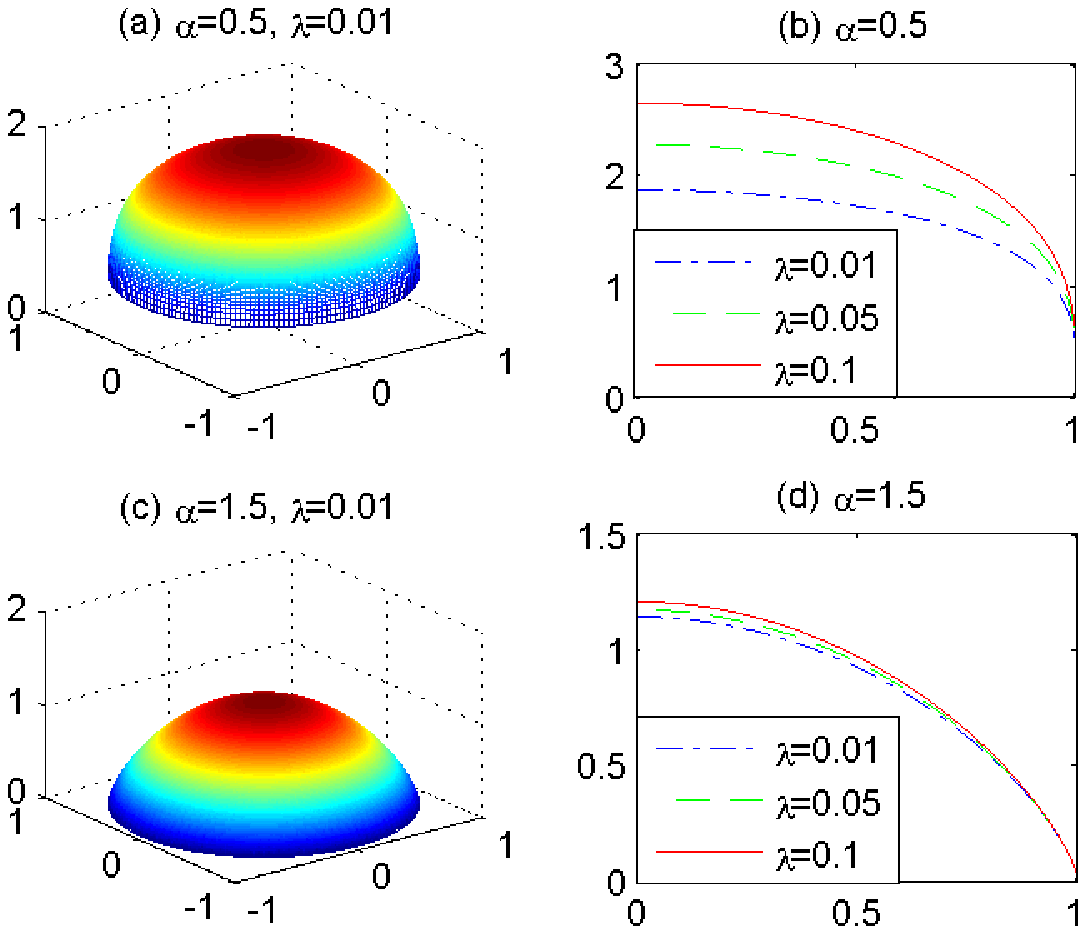}
\caption{MET for the isotropic case  with different $\lam$ and $\a$. }
\label{Met2DsymIS}
\end{center}
\enfig

\par
We use the second-order central differences for $u'(r)$ and $u''(r)$, and take the trapezoidal rule for the nonsingular
 integral terms in (\ref{wholeI}). Assume $f(r)=0$, $d(r)=0$ and $D=B_1(0)$, as the exact solution could not be obtained, we take $U_{640}$ as the 'exact' solution, and $U_J$ is the numerical solution with the resolution $J=640$. Taking
$\lam=0.01, \eps=1$,  we compute the difference between numerical solution $U_J(0)$ and 'exact' solution $U_{640}(0)$
 for $J=10, 20, 40, 80, 160, 320$, i.e., $ error = U_J(0)-U_{640}(0)$ at the fixed point $x=0$. From Fig.~\ref{order2D}, we see that the rate of decay is almost $O(h)$. The Fig.~\ref{Met2DsymIS} (b) and (d) appears the radially symmetric solution of Eq.~(\ref{eq:equivform}). After rotating these two graphs along the vertical axis, we get the mean exit time $u(x,y)$ for $(x,y) \in B_1(0)$ in Fig.~\ref{Met2DsymIS}(a) and (c). When the parameter $\lam$ becomes larger, the `particle' takes more time to exit the domain for these two cases. Moreover, the tempering parameter has more influence for $\a=0.5$ than $\a = 1.5$.

\medskip
\textbf{Acknowledgements}.
 This work was partly supported by the NSFC grant 11901202 (Y.Z.), NSFC grant 11901159 (X.W.),  NSF-DMS no. 1620449 and NSFC grant. 11531006 and
11771449 (J.D.).

\end{document}